\documentclass[12pt]{amsart}
\newtheorem{theorem}{Theorem}[section]
\newtheorem{lemma}[theorem]{Lemma}
\newtheorem{corollary}[theorem]{Corollary}

\newtheorem{remark}[theorem]{Remark}

\newcommand{\ncom}{\newcommand}
\ncom{\lrar}{\longrightarrow}
\ncom{\ov}{\overline}
\ncom{\m}{\mbox}
\ncom{\sta}{\stackrel}
\ncom{\comx}{{\mathbb C}}
\ncom{\A}{{\mathbb A}}
\ncom{\Z}{{\mathbb Z}}
\ncom{\Q}{{\mathbb Q}}
\ncom{\R}{{\mathbb R}}
\ncom{\G}{{\mathbb G}}
\ncom{\al}{\alpha}
\ncom{\p}{{\mathbb P}}
\ncom{\E}{{\mathbb E}}
\ncom{\N}{{\mathbb N}}
\ncom{\K}{{\mathbb K}}
\ncom{\X}{{\mathbb X}}
\ncom{\f}{\frac}
\ncom{\cA}{{\mathcal A}}
\ncom{\cB}{{\mathcal B}}
\ncom{\cD}{{\mathcal D}}
\ncom{\cX}{{\mathcal X}}
\ncom{\cO}{{\mathcal O}}
\ncom{\cW}{{\mathcal W}}
\ncom{\cL}{{\mathcal L}}
\ncom{\cP}{{\mathcal P}}
\ncom{\cH}{{\mathcal H}}
\ncom{\cS}{{\mathcal S}}
\ncom{\cM}{{\mathcal M}}
\ncom{\cC}{{\mathcal C}}
\ncom{\cT}{{\mathcal T}}
\ncom{\cF}{{\mathcal F}}
\ncom{\cN}{{\mathcal N}}
\ncom{\cJ}{{\mathcal J}}
\ncom{\cV}{{\mathcal V}}
\ncom{\cZ}{{\mathcal Z}}
\ncom{\cU}{{\mathcal U}}
\ncom{\cSU}{{\mathcal S \mathcal U}}
\ncom{\cG}{{\mathcal G}}
\ncom{\cQ}{{\mathcal Q}}
\ncom{\cR}{{\mathcal R}}
\ncom{\cY}{{\mathcal Y}}
\ncom{\cE}{{\mathcal E}}

\begin{document}
\baselineskip=16pt
\author[J. N. Iyer]{Jaya NN. Iyer}
\address{The Institute of Mathematical Sciences, CIT
Campus, Taramani, Chennai 600113, India}
\email{jniyer@imsc.res.in}

\title[Chern invariants of some flat bundles ]{Chern invariants of some flat bundles in the arithmetic Deligne cohomology}

\footnotetext{Mathematics Classification Number: 14C25, 14D05, 14D20, 14D21 }
\begin{abstract}
In this note, we investigate the Chern classes of flat bundles in the arithmetic Deligne Cohomology, introduced by Green-Griffiths, Asakura-Saito.
We show nontriviality of the Chern classes in some cases and the proof also indicates that generic flat bundles (in its moduli space) can be expected to have nontrivial classes, which is in contrast with some Gauss-Manin systems.

\end{abstract}
\maketitle

\section{Introduction}

Suppose $X$ is a nonsingular quasi-projective variety defined over the field
of complex numbers. Consider an algebraic vector bundle $\cV$ on $X$ with a connection 
$$\nabla:\cV\lrar \cV\otimes \Omega^1_X$$
 which is flat, i.e. $\nabla^2=0$.
We consider the primary invariants of $(\cV,\nabla)$ 
in various cohomology theories. The Chern classes of $(\cV,\nabla)$ in the Betti, de Rham and the Deligne cohomologies of $X$ are denoted as follows:
$$c^{B}_i(\cV)\,\in\,H^{2i}(X,\Z)$$ 
$$c^{dR}_i(\cV)\,\in\,H^{2i}_{dR}(X,\comx)$$ 
$$c^{D}_i(\cV)\,\in\,H^{2i}_\cD(X,\Z(i)).$$ 

Further, the Deligne Chern class is a lift of the Betti Chern class under the class map
 $$H^{2i}_\cD({X},\Z(i))\lrar H^{2i}({X},\Z).$$ 

By the Chern-Weil theory, the de Rham Chern classes $c^{dR}_i(\cV)$ are zero.
Using the de Rham isomorphism
$$H^\bullet(X,\Z)\otimes \comx \simeq H^\bullet_{dR}(X,\comx)$$
we conclude that the Betti Chern classes $c^{B}_i(\cV)$ are torsion.

When $X$ is projective, it was conjectured by S. Bloch and proved by A. Reznikov (\cite{Re}) that $c^D_i(\cV)$ are torsion, for $i\geq 2$. 
We do not have many examples of flat bundles having nontrivial Chow Chern classes, except 
the case of a direct sum of degree zero line bundles on a general abelian variety (\cite{Bl2}). 

On the other hand, M. Green- P.Griffiths (\cite{Gr-Gri}) and M.Asakura- S. Saito (\cite{As1}) have introduced a category of 
arithmetic Hodge structures. Given a nonsingular variety $X$, they have defined higher Abel-Jacobi maps into some Ext-groups in this category. These groups are termed as the arithmetic Deligne cohomology and denoted by $H^a_{AD}(X,\Q(b))$.
There are cycle class maps

$$CH^p(X)_\Q \sta{\psi^p_{AD}}{\lrar} H^{2p}_{AD}(X,\Q(p))$$

and the usual cycle class map

$$CH^p(X)_\Q \sta{\psi^p_{D}}{\lrar} H^{2p}_{D}(X,\Q(p))$$
factors via $\psi^p_{AD}$ (see also \cite{Sa2}).

The Bloch-Beilinson conjecture on the injectivity of $\psi^p_D$ for varieties over number fields implies the injectivity of the map $\psi^p_{AD}$ for $X/\comx$.

Very few examples are known where non-trivial cycle classes are computed in the arithmetic
Deligne cohomology and which vanish in the Deligne cohomology. Asakura \cite{As2} studied such an explicit zero-cycle on a product of two curves.
 
In this note, we show non-triviality of the Chern classes of some flat bundles in the arithmetic Deligne cohomology.  

We show
\begin{theorem}
Suppose $A$ is a general abelian variety of dimension $g$ defined over $\comx$. Let $E$ be a holomorphic connection on $A$ of rank $n$ and generic in the moduli space of holomorphic connections on $A$.
Then
$$c_i(E)\,\in\,H^{2i}_{AD}(A,\Q(i))$$
is non-zero, for $i\leq \m{min}\{n,g\}$.
\end{theorem}

This is proved in \S 3 Theorem \ref{th.-a}, Remark \ref{re.-hol}.
The proof uses the existence of a flat connection on $E$ and the structure of 
indecomposable flat bundles on $A$, due to Moromito \cite{Mo}. Further, the key point used in the proof is the non-triviality of the Betti Chern classes
of the Poincar\'e bundles.

A similar proof also holds whenever the Betti Chern classes of the Poincar\'e bundle associated to a variety $X$ and a moduli space of flat connections of rank $r$ is nontrivial. This helps to conclude the nontriviality of Chern classes of a generic pair $(X,E)$, $rank(E)=r$ as above.

This gives our second example

\begin{theorem}
Suppose $X=C_{g_1}\times C_{g_2}\times...\times C_{g_n}$ is a product of general smooth
projective curves $C_{g_i}$ of genera $g_i>0$, for $1\leq i\leq n$.
Suppose $E$ is a flat bundle on $X$ and generic in the moduli space of unitary flat 
connections of rank $r$ on $X$. Then
$$c_i(E)\,\in\,H^{2i}_{AD}(A,\Q(i))$$
is non-zero, for $i\leq \m{min}\{n,r\}$.
\end{theorem}

In particular, we obtain the non-triviality of the Chern classes of above flat bundles in the rational Chow groups. Further, these classes lie in the kernel of the Deligne cycle class map, i.e. in the $F^2$- part of the conjectural Bloch-Beilinson filtration.
The non-triviality result is in contrast with the triviality of the Chern classes
of some Gauss-Manin systems (\cite{Mu}, \cite{vdG}, \cite{Iy}, \cite{Es-Vi2}, \cite{Bi-Iy}) in the rational Chow groups, which may be seen as special points in the moduli space of flat connections.

$Acknowledgements:$ We are grateful to U. Jannsen for suggesting to investigate
the Chern classes in the arithmetic Deligne cohomology and for an invitation to visit Regensburg during May 15--June 15, 2005. We also thank U. Jannsen and I. Kausz for having useful discussions. The visit was supported by the DFG project `Algebraic cycle and L-functions'. 

\section{Arithmetic Hodge structures}

We briefly recall the notion of arithmetic Hodge structures and the cycle class map into the arithmetic Deligne cohomology. For further details, we refer to \cite{As1}, \cite{Sa2}.

\subsection{Spreading out a variety $X$}
For any variety $X$ defined over the field of complex numbers, consider a cartesian diagram
\begin{eqnarray*}
X          &\lrar &\cX \\
\downarrow &      &\downarrow \\
\m{Spec}\comx &\lrar & S. 
\end{eqnarray*}
Here $\cX$ and $S$ are defined over a number field and the lower horizontal arrow factors through the generic point of $S$. Such diagrams exist by considering a set of defining equations of $X$ and taking $S$ as the parameter space of the coefficients of these equations. Then $\cX\lrar S$ is called a model of $X$.
M. Saito (\cite{Sa1}) has shown that there is a mixed Hodge module 
$R^i_{dR}(\cX/S)$ on $S$ such that the pullback to Spec $\comx$ is the mixed Hodge structure 
on the cohomology of $X$. The category $MHM(S)$ of mixed Hodge modules over $S$ is an abelian category and having non-trivial higher Ext-groups when $S$ is positive dimensional. 

There is a natural spectral sequence
$$E^{a,b}_1=Ext^b_{MHM(S)}(\Q(c),R^a_{dR}(\cX/S))\Rightarrow Ext^{a+b}_{MHM(\cX)}(\Q(c),\Q).$$    

When $a=2p-k,\,b=k$ and $c=-p$, the Ext-group $Ext^{2p}_{MHM(\cX)}(\Q(-p),\Q)$ 
is the Deligne-Beilinson cohomology of $\cX$. 

The arithmetic Deligne-Beilinson cohomology is defined as the direct limit
of the Ext-groups $Ext^{r}_{MHM(\cX)}(\Q(-s),\Q)$, i.e.,

$$H^r_{AD}(X,\Q(s))\,=\,{lim}_{\cX\lrar S} Ext^r_{MHM(\cX)}(\Q(s),\Q)$$
 and the limit is taken 
over all models of $X$.

Also, notice that 
$$CH^p(X)=lim_{\cX}\, CH^p(\cX)$$ 
since any algebraic cycle on $X$ also spreads out over some $\cX$.

Altogether, the cycle class map into the Deligne-Beilinson cohomology of $\cX$
now yields a cycle class map 
\begin{equation}\label{eq.-2}
CH^p(X)_\Q \sta{\psi^p_{AD}}{\lrar} H^{2p}_{AD}(X,\Q(p))
\end{equation}
and factorising the cycle class map
\begin{equation*}
CH^p(X)_\Q \sta{\psi^p_{D}}{\lrar} H^{2p}_{\cD}(X,\Q(p)).
\end{equation*}

\subsection{Spreading out vector bundles}
 Suppose $E$ is an algebraic vector bundle of rank $r$ on $X$ defined over 
$\comx$.
 Then there exists a finite Zariski open cover $\{U_\al\}$ of $X$ with glueing
 maps 
\begin{equation}\label{eq.-3}
U_\al\cap U_\beta \lrar GL_r(\cO_{U_\al\cap U_\beta}), \,\al\neq \beta.
\end{equation}
 Consider models $\cU_\al\lrar S_\al$, for each $U_\al$ and let
 $\cX=\cup_\al i^{*}_\al\cU_\al\lrar S$ where $S=\cap_\al S_\al$ and $S\sta{i_\al}{\hookrightarrow} S_\al$ is the inclusion. Then $\cX\lrar S$ is a model of $X$. Moreover,
the glueing maps in \eqref{eq.-3} extend to glueing maps
$$\cU_\al\cap \cU_\beta\lrar GL_r(\cO_{\cU_\al\cap \cU_\beta})$$
over an open subset of $S$ and we denote it again by $S$.
In particular, the vector bundle $E$ extends over $\cX\lrar S$ as a family of vector 
bundles $\cE\lrar \cX$. Hence the Deligne Chern class 
$$c^D_i(\cE)\in H^{2i}_\cD(\cX,\Q(i))$$
 extends the Deligne Chern class $c^D_i(E)$.
We call $(\cX\lrar S,\cE)$ as a model for the vector bundle $E$.

\subsection{Arithmetic Deligne Chern class of flat bundles}

Suppose $E$ is an algebraic flat bundle on a nonsingular projective variety $X$
defined over $\comx$.
We denote the image of the Chow Chern classes $c^{CH}_i(E)$ of $E$ under $\psi^i_{AD}$ by
$$c^{AD}_i(E)\,\in \, H^{2i}_{AD}(X,\Q(i)).$$
By functoriality, given any model $\cE$ of $E$, via the natural map $X\lrar \cX$, the Deligne class $c^D_i(\cE)$ pullsback to the arithmetic Deligne class $c^{AD}_i(E)$. 
 Then, by \cite{Re}, the Deligne Chern classes of $E$ are zero. Since the cycle class map $\psi^i_{AD}$ is expected to be injective (\cite{As1}), our aim here
is to show non-triviality of the Chern classes $c^{AD}_i(E)$ in $H^{2i}_{AD}(X,\Q(i))$, for some $(X,E)$.

\begin{remark}\label{re.-1}
Suppose $S'\lrar S$ is a finite surjective cover and 
$f:\cX'=S'\times_S \cX \lrar \cX$ is the pullback morphism of finite degree $d$.
If $\cE$ is a vector bundle on $\cX$ and $\cE'=f^*\cE$ then, by projection formula \cite{Fu},
$$f_*c_i(\cE')\,=\, d.c_i(\cE)$$
in $CH^*(S)_\Q$.
\end{remark}

\section{Holomorphic connections on abelian varieties}

Suppose $E$ is a holomorphic connection of rank $n$ on an abelian variety. Then A. Moromito \cite{Mo} has shown the existence of 
an integrable connection on $E$ (more generally for holomorphic connections on complex torus). This answers positively
a question due to Atiyah. Hence $E$ is a flat bundle on $A$.

A vector bundle on $A$ is said to be \textit{indecomposable} if it cannot be written as a direct sum of two proper subbundles.
         
Let $E=E_1\oplus E_2 \oplus ...\oplus E_r$ be a decomposition of $E$ such that each $E_i$ is an indecomposable vector
bundle admitting a holomorphic connection. Since the Chern character is additive over direct sums, it suffices to compute the
Chern character of indecomposable bundles.

Let $E$ be an indecomposable flat bundle on an abelian variety $A$. Then, by Matsushima \cite{Ma}, we can take
the monodromy representation $\phi$  of $E$ to be 
$$\phi=\sigma\otimes \rho$$
where $\sigma$ is a one dimensional representation and $\rho$ is a unipotent representation. Hence we can write

\begin{equation}\label{eq.-dec}
E=E_\phi=L_\sigma\otimes E_\rho
\end{equation}
where $L_\sigma \,\in\, \m{Pic}^0(A)$ and $E_\rho$ is the flat bundle corresponding to $\rho$.

\begin{lemma}
There is a sequence of vector bundles
$$E_\rho=E_1, E_2,...,E_{l-1},E_l=0$$
such that each $E_i$ has a holomorphic connection and $E_i=E_{i-1}/E^0_{i-1}$. Here $E^0_{i-1}\subset E_{i-1}$ is a 
trivial subbundle, for each $i=1,...,l$.  
\end{lemma}
\begin{proof}
This is Proposition 5.2 and Lemma 5.4 in \cite{Mo}.
\end{proof}

This immediately implies that
\begin{corollary}\label{co.-triv}
In the Grothendieck ring $K(A)_Q$ of $A$,
$$E_\rho=\cO^r$$
for some $r>0$.
\end{corollary}

\begin{corollary}\label{co.-uni}
The Chern character of $E_\phi$ is given as
$$Ch(E_\phi)=r.Ch(L_\sigma) \,\in\, CH^*(A)_\Q$$
\end{corollary}

Hence to compute the Chern character of a holomorphic bundle $E$ on $A$, it suffice to compute the Chern character of
a sum of line bundles of degree zero on $A$, i.e. when $E= \oplus L_i$, for $L_i\in \m{Pic}^0 A$.

We first consider the case when all the $L_i$ are equal to a fixed $L\in \m{Pic}^0(A)$.

\begin{theorem}\label{th.-a}
Suppose $A$ is a generic abelian variety defined over $\comx$. Let $L\in \m{Pic}^0(A)$ be a generic line bundle on $A$. Consider the algebraic flat bundle $E= L^{\oplus n}$ on $A$ and $n\leq g=\m{dim } {A}$.
Then the Chern classes
$$c^{AD}_i(E)\,\in\, H^{2i}_{AD}(A,\Q(p))$$
are non-trivial, for $i\leq n$.
\end{theorem}

In this case, we have
$$c^{CH}_i(E)\,=\, {n\choose i}.(c^{CH}_1(L))^i.$$ 
It was shown in \cite{Bl2} that the classes $(c^{CH}_1(L))^i$ are non-trivial
in $CH^i(A)_\Q$. 

\begin{proof}
Suppose $(\cX\lrar S,\cL)$ is a model for $L$ on $A$. Then $\cX\lrar S$ is an abelian scheme and $\cL\lrar \cX$ is a line
bundle which restricts on any fibre to a degree $0$ line bundle and extending $L\lrar A$. Consider a fine moduli scheme $\cA_g$ of $g$-dimensional polarized abelian varieties (with suitable level structures)
 together with the universal family
$$\A\lrar \cA_g$$
 which are defined over $\Q$.
Consider the family of dual abelian varieties 
$$\m{Pic}^0{\A}\lrar \cA_g$$
and a Poincar\'e line bundle
\begin{eqnarray*}
 & \cP &\\
 &\downarrow & \\
& \A\times_{\cA_g} \m{Pic}^0\A. &\\
\end{eqnarray*}

Furthermore, $\m{Pic}^0\A$ is a fine moduli space defined over $\Q$ of pairs $(A',L')$, where $L'$ is a degree $0$ line bundle on a $g$-dimensional abelian variety $A'$. Hence, by universality, there is a morphism (actually over some finite surjective cover of $S$ which we assume to be $S$, in view of Remark \ref{re.-1})
$$\eta:S\lrar \m{Pic}^0\A $$  
defining the cartesian diagram
\begin{eqnarray*}
\cX &\sta{\eta'}{\lrar} & \A\times_{\cA_g} \m{Pic}^0\A \\
\downarrow\!\beta &  &\,\,\,\,\, \downarrow\!p \\
S &\sta{\eta}{\lrar} &\,\,\,\,\, \m{Pic}^0\A \\
\end{eqnarray*}

 and $\cL={\eta'}^*(\cP)\otimes \beta^*M$.
Also, $\eta$ is defined over the same field of definition as that of $S$ and $M$ is some line bundle on $S$.  

Consider the powers of the first Betti Chern class of $\cP$ 
$$(c^B_1)^i\,\in\,H^{2i}(\A\times_{\cA_g} \m{Pic}^0\A,\Q).$$

Using the degeneration of the Leray spectral sequence for the morphism $\pi:\A\times_{\cA_g} \m{Pic}^0\A\lrar \cA_g$, we have

$$H^{2i}(\A\times_{\cA_g} \m{Pic}^0\A,\Q)= \oplus_{p+q=2i}H^p(\cA_g,R^q\pi_*\Q).$$
Write $(c^B_1)^i= \sum_{p+q=2i} s_{p,q}, \,s_{p,q}\in H^p(\cA_g,R^q\pi_*\Q)$.
 
Since the classes $(c^B_1)^i$ are non-zero on any closed fibre $A'\times \m{Pic}^0A'$ (\cite{Bi-La},p.68-69), the class $s_{0,2i}$ is non-zero, for $i\leq g$.
Hence, $(c^B_1)^i \neq 0$ for $i\leq g$.

Further, using the Leray spectral sequence for the morphism $p$, we have a decomposition

$$H^{2i}(\A\times_{\cA_g} \m{Pic}^0\A,\Q)= \oplus_{p+q=2i}H^p(\m{Pic}^0\A,R^qp_*\Q)$$
and $(c^B_1)^i= \sum_{p+q=2i} s^i_{p,q}, \,s^i_{p,q}\in H^p(\m{Pic}^0\A,R^qp_*\Q)$.

Further, the nondegeneracy of the form of $c^B_i$ (\cite{Bi-La},p.68-69) over any fibre implies that $s^i_{0,2i}$ is a non-zero class, for $i\leq g$.
In particular, there is a non-empty open subset $U\subset \m{Pic}^0\A$ where $s^i_{0,2i}\neq 0$ for each $i$.

Suppose $(A,L)\in U$ is a general point. Then for any model $\cX\lrar S$
of $(A,L)$, the cohomology class $s^i_{0,2i}$, $i\leq g$, which lives on the fibres of $p$ remains non-zero, after base change to $S$. Hence, we now conclude that the class
$$(c^B_1(\cL))^i= ({\eta'}^*c^B_1(\cP)+\beta^*c^B_1(M))^i \in H^{2i}(\cX,\Q)$$ for $i\leq g$, is non-zero, since the powers of the class $\beta^*c^B_1(M)$ do not contribute to the
$(0,2i)$-component.

In particular, the Deligne Chern class 
$$c_i^D(\cL^{\oplus n})= {n\choose i}.(c^D_1(\cL))^i \in H^{2i}_\cD(\cX,\Q(i))$$
is non-zero, for $i\leq n$.
This implies that the arithmetic Deligne Chern class
  $$c^{AD}_i(E)\in H^{2i}_{AD}(A,\Q(i))$$
is non-zero, for $i\leq n$.

\end{proof}

\begin{remark} The same proof also holds if $E$ is a direct sum of distinct bundles $L_i\in \m{Pic}^0(A)$, by considering a larger moduli space parametrizing tuples $(L_1,...,L_n)$, where $L_i\in \m{Pic}^0(A)$. We need to replace the family of Picard groups $\m{Pic}^0(A)$ by the family of self--product ($n$ copies) of Picard groups and repeat the argument.
\end{remark}

\begin{remark}\label{re.-hol}
Using \eqref{eq.-dec}, Corollary \ref{co.-triv} and Corollary \ref{co.-uni}, we deduce that if $E$ is a generic holomorphic 
connection 
on a generic abelian variety then the same conclusion as in Theorem \ref{th.-a} holds for $E$.
\end{remark}

\section{Product of curves}

Let $\cM_g$ denote a fine moduli space of nonsingular projective connected curves of genus $g$, with suitable level structures.

\begin{theorem}
Suppose $C_i$ are general nonsingular projective curves of genus $g_i>0$, for $i=1,...,n$.
Consider the product variety $X=C_{g_1}\times C_{g_2}\times...\times C_{g_n}$ and a flat bundle
$$E=p_1^*E_1\otimes p_2^*E_2 \otimes...\otimes p_n^*E_n$$ on $X$. Here $p_i$ denotes the projection 
to $C_i$ and $E_i$ is a flat bundle on $C_i$ of rank $r_i$.
Suppose $E$ is a generic unitary flat bundle of rank $r=\Pi_{i=1}^n r_i$.
Then the Chern classes of $E$ in the arithmetic Deligne cohomology are non-torsion, i.e.,
$$c^{AD}_i(E)\neq 0 \,\in\,H^{2i}_{AD}(X,\Q(i))$$
for $i\leq \m{min}\{n,r\}$.
\end{theorem}
We denote the rank of $E$ by the $n$-tuple $\bar r=(r_1,...,r_n)$.

\begin{proof}
Consider a fine moduli space of products of curves $C_{g_i}$ of genus $g_i$ (for $1\leq i\leq n$) $$\cM_{\bar g}=\cM_{g_1}\times \cM_{g_2}\times...\times \cM_{g_n}$$ together with the universal product variety
$$\X=\cC_{g_1}\times_{\cM_{\bar g}}...\times_{\cM_{\bar g}} \cC_{g_n}\lrar \cM_{\bar g} .$$

Let $\G_i\lrar \cM_{g_i}$ be the relative moduli space of unitary flat connections of rank $r_i$ and $$\cP_i\lrar \cC_{g_i}\times \G_i$$ be a universal Poincar\'e bundle of rank $r_i$, defined over an \'etale open subset of $\cC_{g_i}\times \G_i$.

Consider the relative moduli space of rank $\bar r$ flat bundles on $\X$ 
\begin{eqnarray*}
&\G=\G_1\times...\times \G_n  & \\
&\downarrow  &\\
&\cM_{\bar g}&
\end{eqnarray*}
together with a universal Poincar\'e bundle defined on an \'etale open subset    
\begin{eqnarray*}
&\cP=p_1^*\cP_1\times\times ...\times p_1^*\cP_n &\\
&\downarrow & \\
& \X \times \cU&\lrar \X \times \G  \\
\end{eqnarray*}

Also all the varieties above are defined over the number field $\ov\Q$.

Suppose $(X,E)$ is any point of $\cU$.
Let $(\cX\lrar S, \cE)$ be a model for $E$ on $X$. Then $\cX\lrar S$ is a family
of products $C_{g_1}\times...\times C_{g_n}$ and $\cE$ restricts on any fibre as a
flat bundle of rank $\bar r$. 
Furthermore, due to universality, there is a morphism 
(actually over a finite surjective cover of $S$ which we assume to be $S$, in 
view of Remark \ref{re.-1})
$$S\sta{\eta}{\lrar} \cU$$
and there is a cartesian diagram (A)

\begin{eqnarray*}
\cX &\sta{\eta'}{\lrar} & \X \times_{\cM_{\bar g}} \cU \\
\downarrow\!\beta &  &\,\,\,\,\, \downarrow\!p \\
S &\sta{\eta}{\lrar} &\,\,\,\,\, \cU \\
\end{eqnarray*}

such that $\cE=M\otimes (\eta ')^*\cP$. 
Here $M$ is the pullback of some line bundle on $S$.

Consider the Betti Chern classes of $\cP$ 
$$c^B_i\,\in\,H^{2i}(\X\times_{\cM_{\bar g}} \cU,\Q).$$

Using the degeneration of the Leray spectral sequence for the morphism $\pi:\X\times_{\cM_{\bar g}} \cU \lrar \cM_{\bar g}$, we 
have

$$H^{2i}(\X\times_{\cM_{\bar g}} \cU,\Q)= \oplus_{k+l=2i}H^k(\cM_{\bar g},R^l\pi_*\Q).$$
Write $c^B_i= \sum_{k+l=2i} s_{k,l}, \,s_{k,l}\in H^k(\cM_{\bar g},R^l\pi_*\Q)$.
 
The classes $c^B_i$ are non-zero on any closed fibre $X\times \cU(X)$, since the 
Kunneth component 
$$ c^B_1(\cP_{j_1})\otimes c^B_1(\cP_{j_2})\otimes...\otimes c^B_1(\cP_{j_i}),\m{ for }j_1<...<j_i$$ of 
$c^B_i(\cP)$ on any fibre is non-zero. Indeed, the class $c^B_1(\cP_{i})\in H^2(C_{g_i}\times \cU(C_{g_i}),\Q)$ is nonzero since the line bundle $\Lambda^{r_i}\cP_i$
is the Poincar\'e line bundle whose class $c_1(\Lambda^{r_i}\cP_i)$ has a nonzero Kunneth component in $H^2(C_{g_i}, \Q)\otimes H^0(\cU(C_{g_i}),\Q)$...(**).
 
We thus deduce that the class $s_{0,2i}$ is non-zero, for 
$i\leq \min\{n,r\}$.
Hence, $c^B_i \neq 0$ for $i\leq \min\{n,r\}$.

Further, using the Leray spectral sequence for the morphism $p$, we have a decomposition

$$H^{2i}(\X\times_{\cM_{\bar g}} \cU,\Q)= \bigoplus_{k+l=2i}H^k(\cU,R^lp_*\Q)$$

and $c^B_i= \sum_{k+l=2i} s^i_{k,l}, \,s^i_{k,l}\in H^k(\cU,R^lp_*\Q)$.

Further, (**) implies that $s^1_{0,2}$ is a non-zero class. Hence $s^i_{0,2i}=(s^1_{0,2})^{\otimes i}$ is a non-zero class, for $i\leq \min\{n,r\}$.

Suppose $(X,E)$ is a general point of $\cU$.
Then the restriction of the non-zero section $s^i_{0,2i}$ to a general point is clearly non-zero. Hence for any model of $(X,E)$ and maps $\eta'$ as in the commutative diagram (A), we now conclude that the class
$$c^B_i(\cE)= {\eta'}^*c^B_i(\cP)+\beta^*c^B_1(M).{\eta'}^*c^B_{i-1}(\cP) \in H^{2i}(\cX,\Q)$$ for $i\leq \min\{n,r\}$, is non-zero, since $c_1^B(M).{\eta'}^*c^B_{i-1}(\cP)$ does not contribute to the $(0,2i)$-Kunneth component.
In particular, the Deligne Chern class 
$$c_i^D(\cE)\,\in\, H^{2i}_\cD(\cX,\Q(i))$$
is non-zero, for $i\leq \min\{n,r\}$.
This implies that the arithmetic Deligne Chern class
  $$c^{AD}_i(E)\in H^{2i}_{AD}(X,\Q(i))$$
is non-zero, for $i\leq \min\{n,r\}$.

\end{proof}

\end{document}